\documentclass{article}
\usepackage{amsmath}
\usepackage{amssymb}

\def\C{{\cal C}}
\def\n{\noindent}
\def\p{\partial}
\def\F{\cal{F}}
 
\def\W{{\cal W}}
\def\Y{{\cal Y}} 
\def\P{\mathbb{P}}
\def\R{\mathbb{R}}

\def\u#1{{\underline{#1}}}

\newtheorem{theo}{Theorem}

\newtheorem{lemma}{Lemma}

\newcommand{\be}{\begin{equation}}
\newcommand{\ee}{\end{equation}}
\newcommand{\ba}{\begin{array}}
\newcommand{\ea}{\end{array}}

\title{Symmetries of special 2-flags}
\author{Piotr Mormul\quad and\quad Fernand Pelletier}
\date{}

\begin{document}
\maketitle
\begin{abstract}
\n This work is a continuation of authors' research interrupted in the year 2010. 
Derived are recursive relations describing for the first time all infinitesimal 
symmetries of special 2-flags (sometimes also misleadingly called `Goursat 2-flags'). 
When algorithmized to the software level, they will give an answer filling in 
the gap in knowledge as of 2010: on one side the local finite classification 
of special 2-flags known in lengths not exceeding four, on the other side 
the existence of a continuous numerical modulus of that classification in 
length seven. 
\end{abstract}
\section{Introduction}\label{Introduction}
The paper is devoted to `special 2-flags', that is, strictly speaking, to rank 3 
distributions generating special 2-flags. More particularly -- to the symmetries of 
such distributions. 
\vskip1mm
There circulates a widely acknowledged folk theorem (cf. section 4 in \cite{Z} and 
p.\,86 in \cite{Mon}) saying that, outside the so-called stable range, distributions 
generic enough do not possess any nontrivial, even only local, symmetry. More to 
the point, in concrete classical classes of subbundles in the tangent bundle, like 
the `3,\,5' or `4,\,7' distributions, the (Lie) groups of symmetries are severely 
restricted in size: not bigger than 14-dimensional in the former (maximal in the 
{\it flat\,} case, when the Cartan tensor -- \cite{C} -- vanishes; \cite{Mon}, p.\,88 
\,and \,\cite{BH}, p.\,456), and not bigger than 21-dimensional in the latter 
(maximal for the {\it instanton\,} distribution, \cite{Mon}, p.\,90). 
And, naturally, likewise restricted in size are the Lie algebras of vector fields 
-- infinitesimal symmetries (i.\,s. for short in all what follows). (They always 
form a Lie algebra due to the Jacobi identity.) While for the objects discussed 
in this work, by virtue of their rather stringent definition, the i.\,s.' algebras 
are {\it infinite-dimensional}. Much like it is the case for the 1-flags, i.\,e., 
{\it Goursat flags} discussed here in length, in the guise of `forerunners', 
in -- still introductory -- Sections \ref{KRw} and \ref{Iso}. (The i.\,s.'s 
for Goursat structures are parametrized by one free function of three variables 
-- a so-called {\it contact hamiltonian}.) 
\vskip1.5mm
The purpose of this paper is to exhibit (for the first time) recursive relations 
which describe all i.\,s.'s of special 2-flags. In order to precise this context 
and give first some motivations, we start from 2-flags of length 1. That is, 
rank 3 distributions $D \subset TM$, $\dim M = 5$ such that $D + [D,D] = TM$ 
(or, the same thing, $[D,D] = TM$; the first order Lie brackets generate all 
the remaining tangent directions; distribution is `two-step'). 
Willy-nilly one enters the domain of the classical `cinq variables' 
work \cite{C}. It was shown there that every such two step $D$ possessed
{\it uniquely\,} determined corank 1 subdistribution $F$ enjoying the property 
\be\label{carr}
[F,\,F] \subset D 
\ee
(see equations (4) on p.\,121 in \cite{C}). Cartan calls such an accompanying 
subdistribution $F$ {\it le syst{\`e}me covariant\,} of \,[the Pfaffian system] $D$. 
Cartan firstly discerns a highly particular situation (a)\,when $[F,F] = F$ 
identically in the vicinity of a point. As a consequence, he infers that, 
in certain local coordinates $t,\,x^0,\,y^0,\,x^1,\,y^1$, $D$ gets description 
$dx^0 - x^1dt = 0 = dy^0 - y^1dt$. In contemporary terminology, such $D$ is, 
up to a local coordinate change, the classical {\it Cartan distribution}, 
or contact system, on the jet space $J^1(1,2)$ of the 1-jets of functions 
$\R(t) \to \R^2(x,\,y)$, with $x^1 = \frac{dx^0}{\!dt}$ and $y^1 = \frac{dy^0}{\!dt}$. 
Its corank 1 covariant subdistribution $F$ (reiterating, involutive in situation (a)!) 
is in these coordinates just ${\rm span}\bigl(\frac{\!\p}{\p x^1}\,,\;\frac{\!\p}{\p y^1}\bigr)$. 
In all what follows we will skip the symbol `span' before a set of vector field generators. 
\vskip1.5mm
By far more interesting is Cartan's situation (b)\,$[F,F] = D$ in the vicinity 
of a given point.\footnote{\,situations (a) and (b) do not exhaust all 
possibilities of the local behaviour of $F$;\\
Elie Cartan used to be interested in clear situations only} 
The covariant object $F$ has then its `curvature' and $D$ is retrievable from 
$F$ alone. We note that situation (b) is extremely rich geometrically and hides 
a {\it functional\,} modulus (one function of five variables) of the local 
classification of `3,\,5' distributions with respect to the diffeos of base 
manifold. 
\vskip1.5mm
We say that a general such $D$ (with no extra information as to (a) or (b)\,) 
generates a 2-flag of length 1, while a $D$ with its covariant system $F$ 
involutive generates a {\it special\,} 2-flag of length 1.\\
Therefore, the adjective `special' in length 1 locally means nothing but 
`jet-like'. How does it look like in bigger lengths/higher jets?
\vskip2mm
Let us analyze the contact system $D$ on a concrete jet space $J^r(1,2) \,=\colon M$ 
with $r \ge 1$. The main observation is that the sequence of modules of vector fields 
-- consecutive Lie squares of $D$, 
\be\label{flag}
TM = D^0 \supset D^1 \supset D^2 \supset \,\cdots \,\supset D^{r-1} \supset D^r, 
\ee
where $D^r = D$ and $[D^j,\,D^j] = D^{j-1}$ for $j = r,\,r-1,\dots,\,2,\,1$, 
grows in ranks regularly by two: $3,\,5,\,7,\dots,\,2r + 1,\,2(r+1) + 1 = \dim M$ 
independently of the underlying points in $M$. (Pay attention to the indexation, which 
starts with the biggest index $r$, following the notation put forward in \cite{MZ1}.) 
The reason is that in passing from $D^j$ to $D^{j-1}$ one forgets about the $j$-th 
order derivatives, so that 
\be\label{rests}
D^{j-1} = \left(D^j,\,\,\frac{\!\p}{\p x^j},\;\frac{\!\p}{\p y^j}\right). 
\ee
Therefore, all these modules of vector fields are actually {\it distributions\,} which 
together form a 2-flag of length $r$ on $M$. Let us scrutinize the members of this flag.\\
The natural coordinates in $J^r(1,2)$ are $t,\,x^0,\,y^0,\,x^1,\,y^1,\dots,\,x^r,\,y^r$, 
where $x^j = \frac{dx^{j-1}}{\!dt}$, $y^j = \frac{dy^{j-1}}{\!dt}$ for $j = 1,\,2,\dots,\,r$. 
In these coordinates the one before last member $D^1$ in \eqref{flag} has a Pfaffian equations' 
description $dx^0 - x^1dt = 0 = dy^0 - y^1dt$, hence it manifestly contains a corank 1 
involutive subdistribution 
\[
F \colon = \,\left(\frac{\!\p}{\p x^j},\;\frac{\!\p}{\p y^j}\,;\ 1 \le j \le r\right). 
\]
Likewise, the next smaller member $D^2$ has description 
\be\label{1.1}
dx^0 - x^1dt = dy^0 - y^1dt = 0 = dx^1 - x^2dt = dy^1 - y^2dt\,,
\ee
hence contains a corank 1 involutive subdistribution 
\[
\left(\frac{\!\p}{\p x^j},\;\frac{\!\p}{\p y^j}\,;\ 2 \le j \le r\right). 
\]
The key point is that the latter happens to be the Cauchy-characteristic module of $D^1$, 
denoted by $L(D^1)$ as in \cite{MZ1}.\footnote{\,\,For $D$ -- a distribution, $L(D)$ is, 
by definition, the module of {\it Cauchy-characteristic} vector fields with values in 
$D$ infinitesimally preserving $D$. That module is automatically (the Jacobi identity) 
closed under the Lie bracket. It is noteworthy that for all the particular distributions 
$D$ occurring in the present work, $L(D) \subset D$ is always not just a module included in $D$, 
but an involutive subdistribution of $D$ of corank 2\,\,(or \,3) when $m = 1\,\,({\rm or}\ \,2)$.}
This pattern replicates itself all the way down the flag. The Pfaffian systems describing $D^j$ 
gradually get larger sets of Pfaffian equations' generators, while the Cauchy-characteristic 
modules get (with a shift in indices!) thinner. In fact, for $1 \le j < r$, 
\[
L(D^j) = \left(\frac{\!\p}{\p x^s},\;\frac{\!\p}{\p y^s}\,;\ j+1 \le s \le r\right)
\] 
sits inside \,$D^{j+1}$ as a corank 1 subdistribution. For instance $L(D^{r-1})$ is 
a field of planes $\bigl(\frac{\!\p}{\p x^r},\;\frac{\!\p}{\p y^r}\bigr)$ sitting 
inside a field of 3-spaces $D^r$, while $L(D^r) = (0)$. Moreover all these geometric 
objects nicely fit together into {\it Sandwich Diagram}, so called after a similar 
(if not identical) diagram assembled for Goursat distributions, or 1-flags, in \cite{MZ1}. 
\[
\!\!\!\!\!\!\!\!\!\!\!\ba{ccccccccccccc}
TM = D^0 & \supset & D^1 & \supset & D^2 & \supset & D^3 & \cdots & 
D^{r-1} & \supset & D^r & & \\
& & \cup & & \cup & & \cup & & \cup & & \cup & & \\
& & F & \supset & L(D^1) & \supset & L(D^2) & \cdots & L(D^{r-2}) & 
\supset & L(D^{r-1}) & \supset & L(D^r) = 0\,. 
\ea
\] 
All vertical inclusions in the diagram are of codimension one, while all 
(drawn, we do not mean superpositions of them) horizontal inclusions are 
of codimension 2. The squares built by these inclusions can, indeed, be 
perceived as certain `sandwiches'. For instance, in the utmost left sandwich 
$F$ and $D^2$ are as if fillings, while $D^1$ and $L(D^1)$ constitute the 
covers (of dimensions differing by 3, one has to admit). At that, the 
sum $2 + 1$ of {\bf co}dimensions, in $D^1$, of $F$ and $D^2$ equals 
the dimension of the quotient space $D^1/L(D^1)$, so that it is natural 
to ask how the 2-dimensional plane $F/L(D^1)$ and the line $D^2/L(D^1)$ 
are mutually positioned in $D^1/L(D^1)$: do they intersect regularly, or 
else the plane subsumes line?\footnote{\,\,this suffices to tell \eqref{1.2} 
from \eqref{1.1} in what follows below} Clearly, that question imposes 
by itself in further sandwiches `indexed' by the upper right vertices 
$D^3,\,D^4,\,\dots,\,D^r$, as well. 
\vskip2mm
\n This question has a trivial answer for the Cartan distribution $D = D^r$ 
analyzed above (all intersections {\it are\,} regular when $r \ge 2$). 
Yet a more pertinent question would be the following. 
\vskip2mm
Assume the existence of Sandwich Diagram with all its above-listed dimensions, 
inclusions, involutivenesses and call such rank 3 distributions $D^r$ \underline{generating 
special 2-flags of length $r$}. Are then those $D^r$ locally `jet-like', that is -- locally 
equivalent to the Cartan contact distribution on $J^r(1,2)$\,?
\vskip2mm
\n For $r = 1$, we reiterate, yes (\cite{C}), but for $r = 2$ already not. 
There suffices to seemingly slightly modify system \eqref{1.1} to 
\be\label{1.2}
dx^0 - x^1dt = dy^0 - y^1dt = 0 = dt - x^2dx^1 = dy^1 - y^2dx^1\,. 
\ee
This rank 3 distribution on $\R^7$ does generate a special 2-flag of length 2, yet is 
{\it not\,} locally equivalent to the `jet-like' one around every point with $x^2 = 0$ 
(cf. \cite{Mor04}, Prop. 1 (iii)). The argument there has been that the object \eqref{1.2} 
has at points $x^2 = 0$ the small growth vector\footnote{\,\,The small growth vector of a 
distribution $D$ at a point $p$ is the sequence of integer numbers $\big(\dim V_j(p)\big)_{j \ge 1}$, 
where $V_1 = D$, \,$V_{j+1} = V_j + [D,\,V_j]$, which {\it ends\,} on the first biggest entry.} 
(3,\,5,\,6,\,7), while the contact system on $J^2(1,2)$ has everywhere the small growth vector 
(3,\,5,\,7). Another, possibly even simpler argument is that at points $x^2 = 0$ there is 
no regular intersection in the only sandwich existing in that length: the line 
$D^2/L(D^1)$ collapses onto the plane $F/L(D^1)$, while the analogous line for 
\eqref{1.1} collapses nowhere.\\
Therefore it follows that the local theory of special multi-flags is not `void' 
in the sense of boling down to the contact systems on the jet spaces for curves. 
In fact, this theory is already fairly rich and still developing, including this work. 
\vskip1mm
Let us reiterate the importance of `special' for 2-flags to be tractable (and the same 
for multi-flags in general). Special, by the way of Sandwich Diagram, brings in so much 
stiffness as to result in the local models with {\bf numerical moduli only}, no functional 
ones. While functional moduli, by simple and widely known dimension counts (cf., for inst., 
section 3 in \cite{Z}) are a commonplace in the local geometry of subbundles in tangent 
bundles. Even the already mentioned paper \cite{C} about 2-flags of length 1 is not yet 
fully understood! On the other side, the initial departing models for us -- contact 
systems on the jet spaces -- are nowadays viewed as just the simplest `baby' 
realizations of the special multi-flags.
\vskip1mm
\n{\it Attention.} This theory is even more neat in that it does not necessitate 
a definition via Sandwich Diagram as such. For it follows from the important works 
\cite{SY,A} that, upon assuming only the properties of the upper row in Sandwich Diagram 
{\it and\,} the existence of a {\it whatever\,} corank one involutive subdistribution 
$F$ in $D^1$, one {\bf automatically} gets Sandwich Diagram in its entirety!\\
In fact, (i)\,such an $F$ is then unique, (ii)\,for $j = 1,\,2,\dots,\,r - 1$ 
there holds 
\[
L(D^j) = D^{j+1} \cap F\,,
\]
(iii)\,$L(D^r) = (0)$ and (iv)\,the $L(D^j)$'s are corank 1 
subdistributions in $D^{j+1}$,\\
so that Sandwich Diagram entirely holds. 
\vskip1mm
Now that the focus is again on Sandwich Diagram, the ongoing question bears 
on the local geometry in the sandwiches `indexed' by the upper right vertices 
$D^2,\,D^3,\,\dots,\,D^r$. It naturally opens the way towards singularities. 
The first step in that direction is a [fairly raw] stratification of germs of 
special 2-flags into so-called sandwich classes. The second is further partitioning 
of sandwich classes into singularity classes (see section \ref{scog}). 
\section{Kumpera-Ruiz watching glasses for Goursat distributions}\label{KRw}
In order to gently introduce the reader to the main techniques of the paper, 
we present in this section a test case -- derive the formulas for the infinitesimal 
symmetries of {\it Goursat\,} distributions which generate 1-flags. This will be 
instrumental during the presentation of similar things to-be-derived for [special] 
2-flags in paper's subsequent sections. 
\vskip1mm
\n Recalling, a rank 2 distribution on a manifold $M$ is Goursat when the tower of 
its consecutive Lie squares, understood as modules of vector fields, consist uniquely 
of regular distributions of ranks \,3,\,\,\,4,\,\,\,5,\,\,\,$\dots$ until $n = \dim M$.
\vskip1mm 
\n With no loss in generality, Goursat distributions understood locally live 
on the stages of Goursat Monster Tower (GMT for short), by some authors called 
alternatively Semple Tower. The stages have been denoted in \cite{MZ2} by 
$\P^r\R^2$, $r \ge 2$. (On the stage $\P^r\R^2$ there lives a Goursat distribution of 
corank $r$.) The best glasses to watch Goursat distributions are Kumpera-Ruiz coordinates 
(KR for short), \cite{KR82}. Those are semi-global sets of coordinates (their domain 
of definition is always dense in a given tower's stage) which critically depend 
on the strata of a most natural stratification of any given stage $\P^r\R^2$ -- 
so-called {\it Kumpera-Ruiz classes}, KR-classes for short, see \cite{MZ1}, p.\,466. 
They exist in $\P^r\R^2$ in number $2^{r-2}$ and are univocally labelled by the words 
of length $r$ over the alphabet $\{1,\,2\}$, with two first letters always 1: 
$1.1.\,i_3.\,i_4.\dots i_r$. (In \cite{MZ1} they were originally labelled by 
the subsets $I \subset \{3,\,4,\dots,\,r\}$, a given $I$ consisting of the indices $j$ 
such that $i_j = 2$.) The KR classes are the main tool in the introductory part 
of our paper. Their generalizations for special 2-flags, so-called singularity 
classes, will play a similar role in the main part of the present contribution 
from Section \ref{resu} onwards. 
\vskip1mm
To each KR-class attached are handy coordinates making that class {\it visible}. 
More precisely, due to the particular topology of the two lowest Monster's stages 
$\P^1\R^2$ and $\P^2\R^2$, they both are unions of pairs of open dense subsets, 
$\P^1\R^2 = U_1 \cup U_2$ and $\P^2\R^2 = V_1 \cup V_2$ such that, for each 
KR-class $\C = 1.1.\,i_3.\,i_4.\dots i_r$ and indices $j,\,k \in \{1,\,2\}$ 
\be\label{UjVk}
\C \cap \pi_{r,1}^{\,-1}(U_j) \cap \pi_{r,2}^{\,-1}(V_k)
\ee
sits in the domain of [Kumpera-Ruiz] coordinates $x_1,\,x_2,\dots,\,x_{r+2}$ 
produced precisely for the data $\C,\,j,\,k$. 
\vskip1.5mm
\n{\bf Remark 1.} The open dense sets $U_j$ and $V_k$ are related to the ways 
the Darboux theorem (in the contact 3D manifold $\P^1\R^2$) and Engel theorem 
(in the {\it Engel\,} 4D manifold $\P^2\R^2$) come into effect. In those coordinates 
\be\label{poly}
\Delta^r = \Big(Y[r],\,\p_{r+2}\Big)\,,
\ee
where, in what follows, $\p_j = \frac{\!\p}{\p x^j}$ and $Y[r]$ is a polynomial 
vector field defined recursively as follows. 

\n Initially $Y[1] = \p_1 + x^3\p_2$ and $Y[2] = Y[1] + x^4\p_3$. When, for $j \ge 3$, 
$Y[j-1]$ is already defined and $i_j = 1$, then $Y[j] = Y[j-1] + x^{j+2}\p_{j+1}$. 
In the opposite case of $i_j = 2$ one puts $Y[j] = x^{j+2}Y[j-1] + \p_{j+1}$. 
The eventual vector field $Y[r]$ in \eqref{poly} is, therefore, polynomial of degree 
(1 + the $\#$ of letters 2 in the code of \,$\C$). That degree is maximal (and equal 
$r - 1$) when the underlying KR-class is 1.1.2.2$\dots$2 ($r - 2$ letters 2 past 
the initial segment 1.1). 
\vskip1.5mm
\n{\bf Remark 2.} Whenever $i_j = 2$ in the code of \,$\C$, the variable $x^{j+2}$ 
brought in at the $j$-th step of the above procedure {\it vanishes\,} at points of 
\eqref{UjVk}. This is a key property of the polynomial visualisations of Goursat 
distributions put forward in \cite{KR82}. 
\vskip1.5mm
The KR-classes are invariant with respect to the local diffeomorphisms of Monster's 
relevant stages. They are only very rough approximations to local models (local normal 
forms). To really approach the orbits, one would need to know the (pseudo-)groups 
of i.\,s.'s of the structures $\Delta^r$ living on $\P^r\R^2$. 
Those groups are infinite-dimensional, for they consist of due prolongations of 
the contact vector fields which preserve the contact structure $\Delta^1$. In order 
to see them, one puts on, no wonder, KR-glasses. That is, works and computes in 
chosen KR-coordinates. 
\section{Infinitesimal symmetries of Goursat flags}\label{Iso}
From now on we assume that KR-coordinates, pertinent for a fixed KR-class in length $r$, 
have been picked and frozen. In these coordinates, every concrete i.\,s. writes down as 
$\Y_f = \sum_{i = 1}^r F^i\p_i$, where the first three components are functions 
of one (smooth) generating function in three variables, say $f(x^1,x^2,x^3)$: 
\be\label{ch}
F^1 = \,-f_3\,,\qquad F^2 = f - x^3f_3\,,\qquad F^3 = f_1 + x^3f_2\,,
\ee
and the remaining components are other, more complicated functions of $f$ 
depending on the KR-class in question, as will be recalled in what follows. 
Such one free function $f$ is called a {\it contact hamiltonian}; the infinite 
dimensionality of the symmetry pseudogroup is visible. 
\vskip1.5mm
When a vector field $\Y_f$ preserves infinitesimally the Goursat $\Delta^r$, 
the {\it truncations\,} of $\Y_f$ do infinitesimally preseve all the earlier 
(older) Goursat structures showing up in the process of building up $\Delta^r$. 
In fact, each component $F^s$, $s = 4,\,5,\dots,\,r+2$, depends {\bf only} 
on the variables $x^1,\,x^2,\dots,\,x^s$ and 
\be\label{j=1.r}
\left[\sum_{i=1}^{j+2}F^i\p_i\,,\ \Delta^j\right] \subset \,\Delta^j
\ee
for $j = 1,\,2,\dots,\,r$, where $\Delta^j = \Big(Y[j],\,\p_{j+2}\Big)$, as in \eqref{poly}. 
This technically central statement is well-known in the theory of Goursat structures, 
compare for instance Proposition 1 in \cite{conham}. Besides, this triangle nature of 
the i.s.'s of Goursat structures will be clearly visible in the recurrences that 
are produced below. The first prolongation of an infinitesimal contactomorphism 
$\sum_{i=1}^3 F^i\p_i$ is $\sum_{i=1}^4 F^i\p_i$, and the new component is 
univocally determined by the previous ones, 
\be\label{F4}
F^4 = Y[2]F^3 - x^4\,Y[2]F^1, 
\ee
compare p.\,222 in \cite{conham}. Reiterating, the components $F^1$ and $F^3$ entering 
formula \eqref{F4} depend on the first three variables, and the field $\overset{\,2}Y$ 
differentiates them accordingly. In the outcome, the component $F^4$ depends on the first 
four variables, and so it goes further on. (This formula is, in fact, subsumed in the line 
of derivations that follow. 
It is given here prior to more involved relations that depend already on the KR-class 
underlying the KR coordinates in use.) 
\vskip1.5mm
We work with a fixed class \,$\C = 1.1.\,i_3.\,i_4.\dots i_r$ and with a fixed letter $i_j$ 
in its code, $j \ge 3$. In order to word the recurrences governing the i.s.'s of $\C$, we need a 
\vskip2mm
\n{\bf Definition} of $s(j)$ for Goursat flags. There can, or cannot, be letters 2 before the letter $i_j$. 
\[
s(j) \colon = \begin{cases} 0\,,&{\rm when\ there\ is\ no\ letter\ 2\ in\ the\ code\ of}\ \C\ {\rm before}\ i_j\,,\\
s\,,&{\rm the\ farthest\ position\ of\ a\ letter\ 2\ before}\ i_j\ {\rm is}\ s,\ {\rm in\ the\ opposite\ case}\,.
\end{cases}
\]
\begin{theo}[\cite{dijon}]\label{the}
Suppose that the components $F^1,\,F^2,\dots,\,F^{j+1}$, $j \ge 3$, of an infinitesimal symmetry 
$\Y_f$ of \,$\Delta^r$ in the vicinty of a KR-class \,$\C = 1.1.i_3.i_4.\dots i_r$ are already 
known. When $i_j = 1$, then 
\[
F^{j+2} = \begin{cases}Y[j] F^{j+1} - \,x^{j+2}Y[2] F^1,&{\rm when}\ s(j) = 0\,,\\
Y[j] F^{j+1} - \,x^{j+2}Y[s(j)] F^{s(j)+1},&{\rm when}\ s(j) \ge 3\,.
\end{cases}
\]
When $i_j = 2$, then 
\[
F^{j+2} = \begin{cases}x^{j+2}\Big(Y[2] F^1 - \,Y[j] F^{j+1}\Big),&{\rm when}\ s(j) = 0\,,\\
x^{j+2}\Big(Y[s(j)] F^{s(j) + 1} - \,Y[j] F^{j+1}\Big),&{\rm when}\ s(j) \ge 3\,.
\end{cases}
\]
\end{theo}
Note before the proof that, on the whole, there are $2^{j-2}$ versions of the formulas for 
the component function $F^{j+2}$, all of them encoded in this theorem. For that many KR-classes 
exist in length $j$. Those formulas are {\it polynomials\,} in the $x$ variables, of growing 
degrees, with coefficients -- partials (of growing orders) of a contact hamiltonian $f$. 
\vskip2mm
The original proof of this theorem occupied full four pages in \cite{dijon}. Now we are going 
to re-prove it in a much shorter manner. Then this new method will be generalized and applied 
to the 2-flags' case in the sections that follow. 
\vskip2mm
\n To begin with, the truncation of the field $\Y_f$ to the Monster level $j$, $\sum_{i=1}^{j+2} F^i\p_i$, 
preserves the Goursat structure $\Delta^j$, as is noted already in \eqref{j=1.r}. Implying, that 
\be\label{ajbj}
\left[\sum_{i=1}^{j+2}F^i\p_i\,,\ Y[j]\right] = \,a_jY[j] + b_j\p_{j+2}
\ee
for certain unspecified functions $a_j$ and $b_j$ of variables $x^1,\dots,\,x^{j+2}$. 
\vskip2mm
\n Now we consider the situation $i_j = 1$. Remembering the construction of 
the field $Y[j]$ when the underlying KR-class is $\C$: 

\n$\bullet$ when $s(j) = 0$, the first ($\p_1$) component on the LHS of \eqref{ajbj} 
is $- Y[2] F^1$. And 

\n$\bullet\bullet$ when $s(j) \ge 3$, the $(s(j) + 1)$-st component on 
the LHS of \eqref{ajbj} is $- Y[s(j)] F^{s(j) + 1}$. So 
\be\label{aj}
a_j = \begin{cases}- Y[2] F^1,&{\rm when}\ s(j) = 0\,,\\
- Y[s(j)] F^{s(j) + 1},&{\rm when}\ s(j) \ge 3\,.
\end{cases}
\ee
One compares now the $(j+1)$-st components on the both sides of \eqref{ajbj}, obtaining 
\[
F^{j+2} - Y[j] F^{j+1} = \,a_j\,x^{j+2}\,.
\]
Substituting on the RHS here the expressions \eqref{aj} in due order, one gets closed 
form formulas for the $\p_{j+2}\,-$\,component function $F^{j+2}$, as invoiced in the theorem. 
As for the coefficient function $b_j$ in \eqref{ajbj}, it is -- here and in what follows 
later -- ascertained last, after finding out $F^{j+2}$. 
\vskip2mm
\n In the situation $i_j = 2$ the arguments differ only technically. Now, regardless 
of the value of $s(j)$, the coefficient $a_j$ can be extracted from \eqref{ajbj} at 
the level $\p_{j+1}$: on the LHS it is $- Y[j] F^{j+1}$, and it is a plain $a_j$ 
on the RHS. Hence 
\be\label{aj'}
a_j = - Y[j] F^{j+1}\,.
\ee
Then, no wonder, one compares the coefficients in \eqref{ajbj} at: $\p_1$, when 
$s(j) = 0$, or else at $\p_{s(j) + 1}$, when $s(j) \ge 3$. In the former case one 
fetches on the LHS the quantity $F^{j+2} - x^{j+2}\,Y[2] F^1$. In the latter, 
the quantity $F^{j+2} - x^{j+2} Y[s(j)] F^{s(j) + 1}$. 
\vskip1.5mm
\n At the same time one fetches $a_j\,x^{j+2}$ on the RHS, just irrelevantly 
of the case in question. That is, accounting for \eqref{aj'}, 
\[
F^{j+2} - x^{j+2}\,Y[2] F^1 = - Y[j] F^{j+1}\,x^{j+2}
\]
(when $s(j) = 0$), or else 
\[
F^{j+2} - x^{j+2}Y[s(j)] F^{s(j) + 1} = - Y[j] F^{j+1}\,x^{j+2}
\]
(when $s(j) \ge 3$). A closed form formula for $F^{j+2}$, invoiced earlier, 
follows immediately. Only then the $b_j$ coefficient is got hold of. In order 
to conclude that the ascertained vector field actually {\it is\,} a symmetry 
of $\Delta^r$ one observes that, in each of the underlying \,$2^{r-2}$ situations, 
\[
\left[\sum_{i=1}^{r+2}F^i\p_i\,,\ \p_{r+2}\right] \,= \,\Big(-\,\p_{r+2}F^{r+2}\Big)\,\p_{r+2}\,,
\]
because only its last component function $F^{r+2}$ depends on the last variable $x^{r+2}$. 
Theorem \ref{the} is now proved. \qquad\qquad$\Box$
\section{Special 2-flags: a basic toolkit}\label{resu}
Special 2-flags constitute a natural follow-up to Goursat flags. 
The latter compactify (in certain precise sense) the contact Cartan distributions 
on the jet spaces $J^r(1,1)$, while the former do the same with respect to 
the jet spaces $J^r(1,2)$.\footnote{\,\,Some researchers, e.g. in \cite{CH}, use, 
instead of `special multi-flags' a somehow misleading synonym `Goursat multi-flags'.}

Sequences of Cartan prolongations of rank 3 distributions are the key players 
in producing (only locally, though) virtually all rank 3 distributions generating 
special 2-flags. There quickly emerges an immense tree of singularities of positive 
codimensions, all of them adjoining the unique open dense Cartan-like strata. 

While the local classification problem is well advanced for the Goursat flags, 
most notably after the work \cite{MZ2}, it is much less advanced for special 2-flags 
(or, more generally, for special multi-flags). It was first attacked in \cite{KR2002}, 
then, in the chronological order, in: \cite{Mor03}, \cite{Mor04}, \cite{V}, \cite{SY}, 
\cite{Mor09}, \cite{A}, and \cite{MP}. After the year 2010 researchers were aiming 
at defining various invariant stratifications in the spaces of germs of special 
multi-flags: \cite{PSl}, \cite{CM}, \cite{CH}, \cite{PS}. The actual state of the art 
is reflected in a recent summarizing work \cite{CCKS}. The works \cite{PSl} and \cite{PS} 
stand out due to a kinematical interpretation of the special 2-flags developed in them. 
Namely, a model of an articulated arm in the 3D space with an engine, or a spacecraft 
with attached string of satellites. The singularities related to various possible 
distributions of {\it right angles\,} between neighbouring segments are already well 
understood and encoded. However, the issue of constructing a kinematics-driven fine 
stratification analogous to Jean's one \cite{J} of the car\,+\,trailers systems 
(modelling 1-flags) in terms of Jean's critical angles, is not yet solved. 
In particular, a faithful expression of the classes in the benchmark work 
\cite{CCKS}, in the terms of an articulated arm in 3D space, seems to be out 
of reach. The issue mentioned above is, most likely, equivalent to that of 
computing all small growth vectors for distributions generating special 2-flags. 

In the work \cite{MP} there was completed only the classification of special 
2-flags in lengths not exceeding 4. At that time the machinery of infinitesimal 
symmetries for those objects was far from being assembled and the techniques 
in use were rather disparate. This notwithstanding, the precise number (34) 
of local equivalence classes of special 2-flags in length 4 was ascertained 
there (cf. the table below). 

The driving force of the present work are the {\it singularity classes} 
(in the occurrence -- of special 2-flags) known for 15 years already. 
They are technically most important for our purposes and results. 
We briefly recall their construction in the next section. For reader's 
convenience, here is the table of cardinalities of singularity classes, 
RV classes of Castro {\it et al\,} \cite{CCKS}, and classes of the local 
equivalence of the special 2-flags, in function of flag's lengths not exceeding 7:  
\begin{center}
\begin{tabular}{|c||c|c||c|}\hline
{\rm length}&$\#$ sing classes&$\#$ RV classes&$\#$ orbits\\ 
\hline
                \hline
$2$ & $2$   & $2$    & $2$  \\ \hline
$3$ & $5$   & $6$    & $7$  \\ \hline
$4$ & $14$  & $23$   & $34$ \\ \hline
\hline
$5$ & $41$  & $98$   & ?    \\ \hline
$6$ & $122$ & $433$  & ??   \\ \hline
$7$ & $365$ & $1935$ & $\infty$ \\ \hline
\end{tabular}
\end{center}
{\it Question.} How to partition a given singularity class of special 2-flags into 
(much finer!) RV classes of \cite{CCKS}\,? \,\,And, all the more so, for special 
$m$-flags, $m > 2$\,?\,!
\subsection{Singularity classes of [germs of] special 2-flags 
refining the sandwich classes}\label{scog}
We first divide all existing germs of special 2-flags of length $r$ 
into $2^{r-1}$ pairwise disjoint {\it sandwich classes\,} in function 
of the geometry of the distinguished spaces in the sandwiches (at the 
reference point for a germ) in Sandwich Diagram on p.\,3, and label 
those aggregates of germs by words of length $r$ over the alphabet 
$\{$1,\,\u{2}$\}$ starting (on the left) with 1, having the second cipher 
\u{2} iff $D^2(p) \subset F(p)$, and for $3 \le j \le r$ having the $j$-th 
cipher \u{2} iff $D^j(p) \subset L(D^{j-2})(p)$. More details about 
the sandwich classes are given in section 1.2 in \cite{MP}. 

\n This construction puts in relief possible non-transverse situations 
in the sandwiches. For instance, the second cipher is \u{2} iff the line 
$D^2(p)/L(D^1)(p)$ is not transverse, in the space $D^1(p)/L(D^1)$, to 
the codimension one subspace $F(p)/L(D^1)(p)$, and similarly in further 
sandwiches. This resembles very much the KR-classes of Goursat germs constructed 
in \cite{MZ1}. \,In length $r$ the number of sandwiches has then been $r-2$ 
(and so the $\#$ of KR classes $2^{r-2}$). For 2-flags the number of sandwiches 
is $r-1$ because the covariant distribution of $D^1$ comes into play and 
gives rise to one additional sandwich. 

Passing to the main construction underlying our present contribution, 
we refine further the singularities of special 2-flags and recall from 
\cite{Mor03} how one passes from the sandwich classes to {\it singularity 
classes}. In fact, to any germ \,$\F$ of a special 2-flag associated is 
a word \,$\W(\F)$ over the alphabet $\{$1,\,2,\,3$\}$, called the `singularity 
class' of \,$\F$. It is a specification of the word `sandwich class' for \,$\F$ 
(this last being over, reiterating, the alphabet $\{$1,\,\u{2}$\}$) with 
the letters \u{2} replaced either by 2 or 3, in function of the geometry 
of $\F$. 
\vskip1.2mm
In the definition that follows we keep fixed the germ of a rank-3 
distribution $D$ at $p \in M$, generating on $M$ a special 2-flag \,$\F$ 
of length $r$. 

\n Suppose that in the sandwich class \,$\C$ of $D$ at $p$ there appears 
somewhere, for the first time when reading from the left to right, the letter 
$\u{2} = j_m$ ($j_m$ is, as we know, not the first letter in \,$\C$) {\bf and} 
that there are in \,$\C$ other letters $\u{2} = j_s$, \,$m < s$, as well. 
We will specify each such $j_s$ to one of the two: 2 or 3. (The specification 
of that first $j_m = \u{2}$ will be made later and will be trivial.) 
Let the nearest \u{2} standing to the left to $j_s$ be $\u{2} = j_t$, 
$m \le t < s$. These two 'neighbouring' letters \u{2} are separated 
in \,$\C$ by $l = s - t - 1 \ge 0$ letters 1.
\vskip1mm
\n The gist of the construction consists in taking the {\it small 
flag\,} of precisely original flag's member $D^s$, 
$$
D^s = V_1 \subset V_2 \subset V_3 \subset V_4 \subset V_5 \subset 
\cdots,
$$
$V_{i+1} = V_i + [D^s,\,V_i]$, then focusing precisely on this new flag's 
member $V_{2l+3}$. Reiterating, in the $t$-th sandwich, there holds the 
inclusion: $F(p) \supset D^2(p)$ when $t = 2$, or else $L(D^{t-2})(p) 
\supset D^t(p)$ when $t > 2$. This serves as a preparation to our 
punch line (cf. \cite{Mor03,Mor09}). 

\n Surprisingly perhaps, specifying $j_s$ to 3 goes via replacing $D^t$ 
by $V_{2l + 3}$ in the relevant sandwich inclusion at the reference point. 
That is to say, $j_s = \u{2}$ is being specified to 3 iff $F(p) \supset 
V_{2l+3}(p)$ (when $t = 2$) or else $L(D^{t-2})(p) \supset V_{2l+3}(p)$ 
(when $t > 2$) holds. 
\vskip1mm
In this way all non-first letters \u{2} in \,$\C$ are, one independently 
of another, specified to 2 or 3. Having that done, one simply replaces 
the first letter \u{2} by 2, and altogether obtains a word over 
$\{1,\,2,\,3\}$. It is the singularity class \,$\W(\F)$ of \,$\F$ at $p$. 
\vskip1.2mm
\n{\bf Example.} In length 4 there exist the following fourteen 
singularity classes: 1.1.1.1, 1.1.1.2; 1.1.2.1, 1.1.2.2, 1.1.2.3; 
1.2.1.1,\footnote{\,\,see section \ref{1211} for more information about 
precisely this class} \,1.2.1.2, \,1.2.1.3, \,1.2.2.1, \,1.2.2.2, \,1.2.2.3, 
\,1.2.3.1, \,1.2.3.2, \,1.2.3.3. (cf. the table on p.\,9).  
\vskip1.2mm
\n(In length $r$ the $\#$ of singularity classes is $\frac12\big(3^{r-1} + 1\big)$; 
the codimension of a class equals the $\#$ of 2's plus twice the $\#$ of 3's 
in the relevant code word.)
\subsection{New approach in the classification problem.}\label{3.2}
A new (2017) approach to the local classification of flags starts with the 
effective (recursive) computation of all infinitesimal symmetries of special 
2-flags, extending the work done (in \cite{dijon}) for 1-flags, reproduced 
with essential shortcuts in Section \ref{Iso} above. The recursive patterns 
depend uniquely on the singularity classes of special 2-flags recapitulated 
above. Those classes are coarser, yes, but much fewer -- see the table preceding 
section \ref{scog} -- than the RV classes summarized (and so neatly systematized) 
in \cite{CCKS}. 
\vskip2mm
Polynomial visualisations of objects in the singularity classes, recalled in 
Section \ref{ekr}, are called EKR's (Extended Kumpera-Ruiz). They `only' feature 
finite families of real parameters. Then the local classification problem is 
rephrased as a search for ultimate normalizations among such families of parameters. 
Having an explicit hold of the infinitesimal symmetries at each prolongation step, 
the freedom in varying those parameters will be ultimately reduced to solvability 
questions of (typically huge) systems of linear equations. 
\vskip1.5mm
In fact, that linear algebra involves only partial derivatives, at the reference 
point, of the first three components of a given infinitesimal symmetry which 
are completely free functions of 3 variables (Lemma \ref{lem0}). 
Keeping the preceding part of a [germ of a] flag in question frozen imposes 
a sizeable set of linear conditions upon those derivatives up to certain order. 
Then some other linear combinations of them appear, or not, to be free -- just 
in function of the local geometry of the prolonged distribution. This, in short, 
would determine the scope of possible normalizations in the new (emerging from 
prolongation) part of EKR's. See sections \ref{1212121} and \ref{1211} below 
for more details. 
\section{EKR glasses for singularity classes of special 2-flags}\label{ekr}
According to section \ref{scog}, the singularity classes of special 2-flags 
of length $r$ are univocally encoded by words of length $r$ over the alphabet 
$\{1,\,2,\,3\}$ such that: -\,the first letter is always 1, and -\,a letter 3, 
if any, must be preceded by a letter 2. 
That is to say, abusing notation a bit, for a singularity class 
\,$\C = 1.i_2.i_3\dots i_r$ over $\{$1,\,2,\,3$\}$, a letter $i_2$ is 
either 1 or 2, and a letter 3 may show up not earlier than at the 3rd position, 
provided there is a letter 2 before it. (We call it, especially in the wider 
context of special $m$-flags with arbitrary $m$, `the least upward jumps rule', 
cf. \cite{Mor04}.)

\n For instance, $\C = 1.2.3$ is a legitimate singularity class of length 3 
(and, in the occurrence, of codimension three in the pertinent Monster's stage No 3). 
\vskip1.5mm
\n For each such \,$\C$ we are going to introduce coordinates, in the number of $2r + 3$, 
\be\label{numer}
t,\,x^0,\,y^0,\,x^1,\,y^1,\dots,\,x^r,\,y^r\,,
\ee
in which the special rank 3 distribution -- let us, from now on, call it $\Delta^r$ 
again -- living on the Monster's $r$-th stage becomes visible. Those coordinates, 
we reiterate it, will sensitively depend on a class \,$\C$. In fact, skipping the 
geometric and also Lie-algebra-related arguments presented in detail in \cite{Mor09}, 
within the domain of those coordinates (subsuming the class \,$\C$), 
\be\label{polly}
\Delta^r = \Big(Z[r],\,\p_{x^r},\,\p_{y^r}\Big)\,,
\ee
where the vector field $Z[r]$ is being defined recursively, shadowing step 
after step the code $1.i_2.i_3\dots i_r$ of \,$\C$. The beginning of 
recurrence is $Z[1] = \p_t + x^1\p_{x^0} + y^1\p_{y^0}$, and, quite simply, 
$\Delta^1 = \Big(Z[1],\,\p_{x^1},\,\p_{y^1}\Big)$ on \,$\R^5(t,x^0,y^0,x^1,y^1)$. 

\n In the recurrence step one assumes description \eqref{polly} known 
for $j - 1$ in the place of $r$, where $1 \le j - 1 \le r - 1$, and puts 
\be\label{trebly}
Z[j] \,= \begin{cases}
Z[j-1] + x^j\p_{x^{j-1}} + y^j\p_{y^{j-1}}\,,&{\rm when}\ i_j = 1\,,\\
x^jZ[j-1] + \p_{x^{j-1}} + y^j\p_{y^{j-1}}\,,&{\rm when}\ i_j = 2\,,\\
x^jZ[j-1] + y^j\p_{x^{j-1}} + \p_{y^{j-1}}\,,&{\rm when}\ i_j = 3\,.
\end{cases}
\ee
In the end of this recurrence (for $j = r$) the description \eqref{polly} 
{\it tout courte} is arrived at, on \,$\R^{2r+3}$ in the variables \eqref{numer}. 
The final first vector field' generator $Z[r]$ is a, possibly deeply involved 
(in function of \,$\C$), polynomial vector field. 
\vskip2mm
Our objective is to ascertain all infinitesimal symmetries $\Y$ of \eqref{polly} 
in the vicinity of any particular class \,$\C$. They will, no wonder, sensitively 
depend on \,$\C$, too. Let us have such $\Y$ expanded in EKR coordinates chosen for \,$\C$:
\be\label{trad}
\Y \,\,= \,A\,\p_t + B\,\p_{x^0} + C\,\p_{y_0} + \sum_{s = 1}^r \Big(F^s\,\p_{x^s} + G^s\,\p_{y^s}\Big)\,.
\ee
The first key property (needed later) is 
\begin{lemma}\label{lem0}
The component functions $A,\,B,\,C$ in {\rm \eqref{trad}} 
depend only on the variables $t,\,x^0,\,y^0$.
\end{lemma}
Proof of Lemma \ref{lem0}. The reason is that, whatever the class \,$\C$, 
in the chosen EKR coordinates associated to \,$\C$ the bottom row in Sandwich 
Diagram has formally the same description as for the Cartan contact system 
on $J^r(1,2)$. In particular, because the relations \eqref{rests} keep 
holding true in the vicinity of \,$\C$ in these coordinates, the covariant 
subdistribution $F$ of $D^1$ is there invariably of the form 
\[
F \,= \,\left(\p_{x^i},\,\,\p_{y^i}\,;\ 1 \le i \le r\right)
\]
The symmetry $\Y$, preserving $\Delta^r = \colon D$, preserves the derived flag 
\,$\big(D^j\big)_{j = r}^{\,\,0}$ \,of \,$D$, so preserves this $F$, too. Hence 
the first three components of $\Y$ can{\it not\,} depend on the variables $x^i$ 
and $y^i$ for $1 \le i \le r$, as stated in the lemma.\qquad\quad$\Box$
\vskip2mm
\n{\bf Remark 3.} Note, however, one essential difference with the 1-flags in 
that here are {\bf three} free functions in the base of the theory, instead 
of just one contact hamiltonian there (in formulas \eqref{ch}\,). 
\vskip2mm
As previously, one needs some additional information about the code of \,$\C$. 
So for $j = 2,\,3,\dots,\,r$ we define 
\[
s(j) = \begin{cases} 0\,,&{\rm when}\ i_2,\dots,\,i_{j-1} = 1\,,\\
\max\{s \colon \,2 \le s < j\ \,\&\ \,i_s > 1\}\,,&\textrm{in the opposite case}\,.
\end{cases}
\]
Note that when $s(j) \ge 2$, then $i_{s(j)} = 2$ or else $i_{s(j)} = 3$. 
These two distinct (and disjoint) geometric situations account for bigger complexity 
of the recurrences to be produced. (The eventail of possible singularities of 
special 2-flags is much wider than for Goursat.)
\section{Infinitesimal symmetries of special 2-flags got hold of}\label{got}
Our main theorem of the paper, Theorem \ref{THE} below, shows that every infinitesimal 
symmetry is uniquely determined by the singularity class under consideration together 
with symmetry's first three component functions, denoted traditionally $A,\,B,\,C$, 
{\bf in an explicit, algorithmically computable manner}. Namely, 
\begin{theo}\label{THE}
Let $U$ be the domain of EKR coordinates \eqref{numer} chosen for an arbitrarily fixed 
singularity class $1.\,i_2.\,i_3\dots.\,i_r$. In those coordinates, all infinitesimal symmetries 
\,$\Y$ of \,$\Delta^r$ restricted to $U$ are of a particular form \eqref{trad}, where $A,\,B,\,C$ 
are free smooth functions of only $t,\,x^0,\,y^0$ and the $F^s,\,G^s$, $1 \le s \le r$, are 
univocally recursively determined by $A,\,B,\,C$ and the class code, according to the formulae 
given in \eqref{F1G1} and Lemmas \,{\rm \ref{lem1}, \,\ref{lem2}} and \,{\rm \ref{lem3}} below. 
\end{theo}
PROOF. We are going to ascertain one by one (or rather two by two) the consecutive components 
of vector fields $\Y$ in \eqref{trad} above, from $F^1$ and $G^1$ on, given the initial 
arbitrary function data $A,\,B,\,C$. To this end we will use the truncations $\Y[j]$ of 
$\Y$ to the spaces of coordinates of indices $\le j$, $j = 1,\,2,\dots,\,r$, on which 
the distributions $\Delta^j$ live:  
\be\label{trun}
\Y[j] \,= \,A\,\p_t + B\,\p_{x^0} + C\,\p_{y_0} + \sum_{s = 1}^j \Big(F^s\,\p_{x^s} + G^s\,\p_{y^s}\Big).
\ee
{\it Attention}. The formulas \eqref{F1G1} right below and in Lemmas \ref{lem1}, \ref{lem2} and 
\ref{lem3} below are, in the first place, only {\it necessary\,} for $\Y$ to be a true symmetry 
of $\Delta^r$. They became also {\it sufficient\,} in the last part of our (long) proof of 
Theorem \ref{THE}. 
\vskip2mm
To begin with, let us demonstrate the argument on the `baby' components $F^1$ and $G^1$. 
The infinitesimal invariance condition 
\[
\left[\Y[1]\,,\ \Delta^1\right] \subset \Delta^1
\]
clearly implies 
\be\label{baby}
\left[\Y[1]\,,\ Z[1]\right] = a_1\,Z[1] + b_1\p_{x^1} + c_1\p_{y^1}\,,
\ee
which in turn implies $a_1 = - Z[1] A$. At the same time $F^1 - Z B[1] = a_1\,x^1$ 
and $G^1 - Z[1] C = a_1\,y^1$. Putting all this together, 
\be\label{F1G1}
\begin{cases} F^1 = \,Z[1] B - x^1Z[1] A\,, & \\
G^1 = \,Z[1] C - y^1Z[1] A\,. &
\end{cases}
\ee
So indeed the pair of new components in $\Y[1]$ is univocally determined by the base 
components $A,\,B,\,C$. As for the coefficients $b_1$ and $c_1$ in \eqref{baby}, 
they get ascertained only {\it after\,} $F^1$ and $G^1$ are found.\\
This inference is an instance of a general 
\begin{lemma}\label{lem1}
Assuming that an infinitesimal symmetry $\Y[j-1]$ of \,$\Delta^{j-1}$ is already known 
for certain $2 \le j \le r$, {\bf in the situation} $i_j = 1$, the $\p_{x^j}\,-$ and 
$\p_{y^j}\,-$\,components of the prolongation $\Y[j]$ of \,$\Y[j-1]$ are as follows 
\[
F^j \,= \begin{cases} Z[j] F^{j-1} - \,x^j Z[1] A\,,&{\rm when}\ s(j) = 0\,,\\
Z[j] F^{j-1} - \,x^j Z[s(j)] F^{s(j) - 1}\,,&{\rm when}\ s(j) \ge 2,\ i_{s(j)} = 2\,,\\
Z[j] F^{j-1} - \,x^j Z[s(j)] G^{s(j) - 1}\,,&{\rm when}\ s(j) \ge 2,\ i_{s(j)} = 3\,.
\end{cases}
\]
\[
G^j \,= \begin{cases} Z[j] G^{j-1} - \,y^j Z[1] A\,,&{\rm when}\ s(j) = 0\,,\\
Z[j] G^{j-1} - \,y^j Z[s(j)] F^{s(j) - 1}\,,&{\rm when}\ s(j) \ge 2,\ i_{s(j)} = 2\,,\\
Z[j] G^{j-1} - \,y^j Z[s(j)] G^{s(j) - 1}\,,&{\rm when}\ s(j) \ge 2,\ i_{s(j)} = 3\,.
\end{cases}
\]
\end{lemma}
Proof of Lemma \ref{lem1}. The vector field $\Y[j]$ infinitesimally preserves 
the distribution $\Delta^j$, whence 
\be\label{arj}
\left[\Y[j]\,,\;Z[j]\right] = a_j Z[j] + b_j\,\p_{x^j} + c_j\,\p_{y^j}
\ee
for certain unspecified functions $a_j,\,b_j,\,c_j$. The coefficient $a_j$ is of central importance 
here. We typically work, here and in what will follow later, in the following order: -\,we firstly 
ascertain $a_j$, -\,secondly find (this is most important) $F^j$ and $G^j$, -\,eventually ascertain 
the values of $b_j$ and $c_j$. 

\n The function $a_j$ can be extracted from \eqref{arj} by watching this vector equation on the level 
of such a component of $Z[j]$ which is identically 1. Inspecting the stepwise construction that 
leads from $Z[1]$ to $Z[j]$, there always {\it is\,} such a component! Namely, it is the 
$\p_t\,-$\,component when $s(j) = 0$. When, on the contrary, $s(j) \ge 2$, it is either 
the $\p_{x^{s(j) - 1}}\,-$\,component (when $i_{s(j)} = 2$), or else it is the 
$\p_{y^{s(j) - 1}}\,-$\,component (when $i_{s(j)} = 3$). With thus specified 
information, it is a matter of course that 
\be\label{ajj}
a_j \,= -\begin{cases} Z[1] A\,,&{\rm when}\ s(j) = 0\,,\\
Z[s(j)] F^{s(j) - 1}\,,&{\rm when}\ s(j) \ge 2,\ i_{s(j)} = 2\,,\\
Z[s(j)] G^{s(j) - 1}\,,&{\rm when}\ s(j) \ge 2,\ i_{s(j)} = 3\,.
\end{cases}
\ee
On the other hand, the same equation \eqref{arj} watched on the level of $\p_{x^{j-1}}$ reads 
\[
F^j - Z[j] F^{j-1} = \,a_j\,x^j\,,
\]
and watched on the level of $\p_{y^{j-1}}$ reads 
\[
G^j - Z[j] G^{j-1} = \,a_j\,y^j\,. 
\]
The needed expressions for $F^j$ and $G^j$ follow upon substituting the expression 
\eqref{ajj} of $a_j$ into these two equations.\qquad\qquad$\Box$
\begin{lemma}\label{lem2}
Assuming that an infinitesimal symmetry $\Y[j-1]$ of \,$\Delta^{j-1}$ is already known 
for certain $2 \le j \le r$, {\bf in the situation} $i_j = 2$, the $\p_{x^j}\,-$ and 
$\p_{y^j}\,-$\,components of the prolongation $\Y[j]$ of \,$\Y[j-1]$ are as follows 
\[
F^j \,= \begin{cases}x^j\left(Z[1] A - Z[j] F^{j-1}\right)\,,&{\rm when}\ s(j) = 0\,,\\
x^j\left(Z[s(j)] F^{s(j) - 1} - Z[j] F^{j-1}\right)\,,&{\rm when}\ s(j) \ge 2,\ i_{s(j)} = 2\,,\\
x^j\left(Z[s(j)] G^{s(j) - 1} - Z[j] F^{j-1}\right)\,,&{\rm when}\ s(j) \ge 2,\ i_{s(j)} = 3\,.
\end{cases}
\]
\[
G^j = Z[j] G^{j-1} - \,y^j Z[j] F^{j-1}\,. 
\]
\end{lemma}
Proof of Lemma \ref{lem2}. The vector equation \eqref{arj} still holds true. 
Now the $a_j$ coefficient can be (and easily) extracted from it at the level $\p_{x^{j-1}}$, 
because the coefficient of the $\p_{x^{j-1}}\,-$\,component in $Z[j]$ is 1: 
\be\label{ajjj}
a_j \,= - Z[j] F^{j-1}. 
\ee
At the same time writing down the equal sides of \eqref{arj} at the level $\p_{y^{j-1}}$, 
\[
G^j - Z[j] G^{j-1} = \,a_j\,y^j, 
\]
leads, by the way of \eqref{ajjj}, to the desired formula for $G^j$. 
\vskip1.5mm
\n It is not that quick with the function $F^j$. It can be extracted from precisely one 
out of three levels of the $\p_t\,-$, $\p_{x^{s(j)-1}}\,-$, or $\p_{y^{s(j)-1}}\,-$\,components. 
Because one, once again, looks for a component in $Z[j]$ with a coefficient 1, if `enveloped' 
now in the factor $x^j$ (because $i_j > 1$ in the proposition under proof).  

\n In function of the position of that `1', equalling the relevant levels in \eqref{arj}, 
one gets precisely one relation out of the following three 
\[
\begin{cases}F^j - x^j Z[1] A = a_j\,x^j\,,&{\rm when}\ s(j) = 0\,,\\
F^j - x^j Z[s(j)] F^{s(j) - 1} = a_j\,x^j\,,&{\rm when}\ s(j) \ge 2,\ i_{s(j)} = 2\,,\\
F^j - x^j Z[s(j)] G^{s(j) - 1} = a_j\,x^j\,,&{\rm when}\ s(j) \ge 2,\ i_{s(j)} = 3\,.
\end{cases}
\]
Then, accounting for \eqref{ajjj}, the desired formula for $F^j$ follows.\qquad\qquad$\Box$
\begin{lemma}\label{lem3}
Assuming that an infinitesimal symmetry $\Y[j-1]$ of \,$\Delta^{j-1}$ is already known 
for certain $2 \le j \le r$, {\bf in the situation} $i_j = 3$, the $\p_{x^j}\,-$ and 
$\p_{y^j}\,-$\,components of the prolongation $\Y[j]$ of \,$\Y[j-1]$ are as follows 
\[
F^j \,= \begin{cases} x^j\left(Z[1] A - Z[j] G^{j-1}\right)\,,&{\rm when}\ s(j) = 0\,,\\
x^j\left(Z[s(j)] F^{s(j) - 1} - Z[j] G^{j-1}\right)\,,&{\rm when}\ s(j) \ge 2,\ i_{s(j)} = 2\,,\\
x^j\left(Z[s(j)] G^{s(j) - 1} - Z[j] G^{j-1}\right)\,,&{\rm when}\ s(j) \ge 2,\ i_{s(j)} = 3\,.
\end{cases}
\]
\[
G^j = Z[j] F^{j-1} - \,y^j Z[j] G^{j-1}\,.
\]
\end{lemma}
Proof of Lemma \ref{lem3}. Invariably, the vector equation \eqref{arj} keeps holding true. 
The $a_j$ coefficient on its right hand side can be extracted from it at the level $\p_{y^{j-1}}$, 
because now the coefficient of the $\p_{y^{j-1}}\,-$\,component in $Z[j]$ is 1: 
\be\label{ajjjj}
a_j \,= - Z[j] G^{j-1}. 
\ee
Then, writing simply down the equal sides of \eqref{arj} at the level $\p_{x^{j-1}}$, 
\[
G^j - Z[j] F^{j-1} = \,a_j\,y^j, 
\]
leads, by the way of \eqref{ajjjj}, to the presently needed formula for $G^j$. 
\vskip1.5mm
\n As for the function $F^j$, it can again be extracted from precisely one out of 
three levels of the $\p_t\,-$, $\p_{x^{s(j) - 1}}\,-$, or $\p_{y^{s(j) - 1}}\,-$\,components. 
In function of the position of that key component `1' in the field $Z[j]$, equalling 
the sides of the relevant levels in \eqref{arj}, one gets precisely one relation 
out of the following three 
\[
\begin{cases}F^j - \,x^j Z[1] A \,= \,a_j\,x^j\,,&{\rm when}\ s(j) = 0\,,\\
F^j - \,x^j Z[s(j)] F^{s(j) - 1} \,= \,a_j\,x^j\,,&{\rm when}\ s(j) \ge 2,\ i_{s(j)} = 2\,,\\
F^j - \,x^j Z[s(j)] G^{s(j) - 1} \,= \,a_j\,x^j\,,&{\rm when}\ s(j) \ge 2,\ i_{s(j)} = 3\,.
\end{cases}
\]
Upon accounting for \eqref{ajjjj}, the expected formula for $F^j$ follows.
\qquad\qquad$\Box$
\vskip2mm
As already invoiced, the obtained recursive formulas -- at this moment only necessary 
-- are also sufficient for the produced vector field $\Y$ to actually {\it be} a symmetry 
of $\Delta^r$. Indeed, knowing already that \,$[\Y,\,Z[r]] \in \Delta^r$ (cf. the always 
holding true formulas \eqref{arj} taken now for $j = r$), what only remains to be done 
is to take the remaining two generators of $\Delta^r$ and justify the vector fields' 
inclusions 
\[
\big[\Y,\,\p_{x^r}\big],\ \,\big[\Y,\,\p_{y^r}\big] \,\in \,\Delta^r. 
\]
To that end we note that Lemma \ref{lem0} coupled with formulas \eqref{F1G1} and 
all those listed in auxiliary Lemmas \ref{lem1}, \ref{lem2} and \ref{lem3} 
yield by simple induction that, for $j = 1,\,2,\dots,\,r$, 
\[
\textrm{the components $F^j$ and \,$G^j$ of \,$\Y$ depend only on $t,\,x^0,\,y^0,\,x^1,\,y^1,\dots,\,x^j,\,y^j$}. 
\]
Using this information for $1 \le j \le r - 1$ and again Lemma \ref{lem0}, 
one computes with ease 
\[
\big[\Y,\,\p_{x^r}\big] \,\,= \,-\,\big[\p_{x^r},\,\Y\big] \,= \big(-\,\p_{x^r}\,F^r\big)\,\p_{x^r} + 
\big(-\,\p_{x^r}\,G^r\big)\,\p_{y^r}
\]
and 
\[
\big[\Y,\,\p_{y^r}\big] \,\,= \,-\,\big[\p_{y^r},\,\Y\big] \,= \big(-\,\p_{y^r}\,F^r\big)\,\p_{x^r} + 
\big(-\,\p_{y^r}\,G^r\big)\,\p_{y^r}.
\]
Now, at long last, the proof of Theorem \ref{THE} is complete.\qquad\qquad$\Box$
\section{Applications of recursively computable infinitesimal symmetries 
to the local classification problem}
The main motivation underlying the present contribution has been to advance results 
in the local classification problem for special 2-flags -- to propose a late follow-up 
to the work \cite{MP}. In fact, getting -- recursively -- hold of the infinitesimal 
symmetries of special 2-flags\footnote{\,and, as a matter of fact, of all special 
$m$-flags, $m \ge 2$, too -- this being the subject of a possible another paper} opens 
a way to advance the local classification in lengths $r = 5$ (cf. in this respect, 
in particular, section \ref{1211}) and $r = 6$ which have kept challenging the small 
monster community for the last 15 years (see the table preceding section \ref{scog}). 
\subsection{Continuous modulus in the class 1.2.1.2.1.2.1}\label{1212121}
\n Reiterating already, the exact local classification of special 2-flags (and, all 
the more so, all special {\it multi}-flags) in lengths exceeding 4 is, in its generality, 
unknown. It is not excluded that a continuous modulus of the local classification hides 
itself already somewhere in length 6. Instead, we want to give an example in length 7 
of the effectiveness of our formulae put forward in Section \ref{got}. 
\vskip1.5mm
\n A possibly deepest fact communicated in \cite{MP} was 
\begin{theo}[\cite{MP}]\label{l7}
In the singularity class \,$\C = 1.2.1.2.1.2.1$ of special \,$2$-flags of 
length $7$ there resides a continuous modulus of the local classification. 
\end{theo}
This was originally proved (in the year 2003, as a matter of fact) by brute 
force, and here is how the i.s.'s may help. 
\vskip2mm
\n PROOF. In the coordinates constructed for the class \,$\C$ we work with certain 
{\it germs\,} of the distribution $\Delta^7$ which generates a locally universal 
special 2-flag of length 7. The reference points for those germs  belong to \,$\C$. 
More precisely, these are the points, say $P$, with the coordinates 
\begin{eqnarray}\label{point} 
&t = x^0 = y^0 = x^1 = y^1 = x^2 = y^2 = 0\,,\quad x^3 = 1\,,\\
&y^3 = x^4 = y^4 = 0\,,\ \ x^5 = 1\,,\ \ y^5 = x^6 = y^6 = 0\,,\ \ x^7 = c\,,\ y^7 = 0\,.\nonumber
\end{eqnarray}
We {\it intend\,} to infinitesimally move such $P$ only in the $\p_{x^7}-$\,direction. 
(Compare, for instance, \cite{MZ1}, where also only the farthest part of a flag -- Goursat 
in that occurrence -- was subject to possible movies.) That is, we look for an i.\,s. 
having at a point $P$ of type \eqref{point} {\bf all but the $\p_{x^7}-$\,components 
vanishing}. Remembering about the triangle pattern of dependence of those component functions, 
this means the vanishing of $A,\,B,\,C$ at $(0,0,0)$, the vanishing of $F^j\big(\pi_{7,j}(P)\big),\,
G^j\big(\pi_{7,j}(P)\big)$ for $j = 1,\,2,\dots,\,6$ and the vanishing of $G^7(P)$. 
The component $F^7(P)$ is not yet known and will be analyzed with care. 

\n Initially we do not know how few/many such vector fields could exist. At any rate, 
any one of them is induced by certain functions $A,\,B,\,C$ in the variables $t,\,x,\,y$. 
The recurrence formulae are known from Section \ref{got}. When, among other components 
of an i.\,s., one wants to express $F^7(P)$ via those basic unknown functions $A,\,B,\,C$, 
one goes {\it backwards\,} along the code of \,$\C$, and firstly applies Lemma \ref{lem1} 
(because $i_7 = 1$), then Lemma \ref{lem2} (because $i_6 = 2$), then again Lemma \ref{lem1} 
(because $i_5 = 1$), and so on intermittently. Upon applying with care these lemmas due 
numbers of times, the above-listed vanishings mean in the terms of the functions 
in the base 
\begin{eqnarray*}
0 \,= A(0,0,0) = B(0,0,0) = C(0,0,0) = &B_t(0,0,0) = C_t(0,0,0)\\
&= C_{x^0}(0,0,0) = c\,C_{t\,x^0}(0,0,0)\,,
\end{eqnarray*}
and -- most important 
\be\label{F3}
0 \,= F^3\big(\pi_{7,3}(P)\big) = \big(3A_t - 2B_{x^0}\big)(0,0,0)\,,
\ee
\be\label{F5}
0 \,= F^5\big(\pi_{7,5}(P)\big) = \big(B_{x^0} - A_t\big)(0,0,0)\,. 
\ee
Now comes the punch line, because the outcome of the computations for $F^7$ is 
\be\label{F7}
F^7(P) = 3c\big(A_t - B_{x^0}\big)(0,0,0)\,. 
\ee
Relations \eqref{F3} and \eqref{F5} together imply $A_t(0,0,0) = B_{x^0}(0,0,0) = 0$. 
So $F^7(P) = 0$ by \eqref{F7}. That is, every i.\,s. of \,$\C$ must infinitesimally freeze 
at $P$ the coordinate $x^7$, when it infinitesimally freezes all the remaining coordinates 
specified in \eqref{point}. Theorem \ref{l7} is proved.\qquad\qquad$\Box$ 
\vskip2mm
\n{\bf Remark 4.} In other terms, the germs of the structure $\Delta^7$ at various 
points $P$ as above (i.\,e., for different values of the parameter $c$) are pairwise 
non-equivalent. The local geometry of the distribution $\Delta^7$ changes continuously 
within the discussed class \,$\C$. 
\subsection{Towards the classification of the one step 
prolongations within singularity class 1.2.1.1}\label{1211}
We conclude the paper by excerpting from \cite{MP} the partition, into the orbits of 
the local classification, of the singularity class 1.2.1.1 (when the width $m = 2$, 
cf. Remark 5 on p.\,37 there), and suggesting a line of possible continuation in 
the next length 5. This class is not chosen at random; it splits into maximal (6) 
number of orbits in that length 4, cf. Section 7 in \cite{MP}. The names of orbits 
are taken from that preprint. One means {\bf the germs of} $\Delta^4$, watched in 
the EKR coordinates constructed for 1.2.1.1, {\bf at points}, say $P$, having 
$t = x^0 = y^0 = x^1 = y^1 = x^2 = y^2 = 0$ \,and 
\begin{center}
\begin{tabular}{|c|||c|c||c|c|}\hline
{\rm the orbit}&$x^3\Big(\pi_{4,3}(P)\Big)$&$y^3\Big(\pi_{4,3}(P)\Big)$&$x^4(P)$&$y^4(P)$\\ 
\hline
                \hline
$1.2.1_{-\rm s,tra}.1$              & $1$ & $0$ & $0$ & $0$ \\ \hline
$1.2.1_{-\rm s,tan}.1_{-\rm s,tra}$ & $0$ & $1$ & $1$ & $0$ \\ \hline
$1.2.1_{-\rm s,tan}.1_{-\rm s,tan}$ & $0$ & $1$ & $0$ & $0$ \\ \hline
$1.2.1_{+\rm s}.1_{-\rm s,tra}$     & $0$ & $0$ & $1$ & $0$ \\ \hline
$1.2.1_{+\rm s}.1_{-\rm s,tan}$     & $0$ & $0$ & $0$ & $1$ \\ \hline
$1.2.1.1_{+\rm s}$                  & $0$ & $0$ & $0$ & $0$ \\ \hline
\end{tabular}
\end{center}
\vskip.1cm
Upon prolonging $\Delta^4$ to $\Delta^5$ in the vicinity of points of 1.2.1.1, 
one is to work with points in the classes $1.2.1.1.i_5$, $i_5 \in \{1,\,2,\,3\}$. 
The classification result recalled in the table above applies now to the distribution 
$\left[\Delta^5,\,\Delta^5\right]$ and as such remains true, regardless of the value 
of $i_5$ (the Lie square of $\Delta^5$ does {\it not\,} depend on new variables $x^5,\,y^5$). 
The same concerns the recursive formulae for the component functions $F^j,\,G^j$, $j = 1,2,3,4$, 
of the i.s.'s of $\Delta^5$. Yet, naturally, expressions for the components $F^5,\,G^5$ 
depend critically on the value of $i_5$. Sticking to the points $P$ from the table, 
one is to analyze the expressions for $F^5(Q)$ and $G^5(Q)$, $Q \in 1.2.1.1.i_5$, 
$\pi_{5,4}(Q) = P$. They are linear in $x^5(Q),\,y^5(Q)$, with coefficients depending 
on $P$ and on certain partials at $(0,0,0)$ of the basic functions $A,\,B,\,C$. 
All the difficulty resides in the -- unknown and hard to compute -- coefficients 
standing next to those partials. 
\vskip1.5mm
\n An instructive example is given in section \ref{1212121}. The coefficient 
standing next to $c = x^7(P)$ on the RHS of \eqref{F7} has appeared {\it forced\,} 
to be zero by the earlier infinitesimal normalizations \eqref{F3} and \eqref{F5}. 
Because of that phenomenon, even the outcome of the classification of singularity 
class 1.2.1.1.1 ($i_5 = 1$) is difficult to predict. 
\vskip1.5mm
In general -- in higher lengths -- {\it systems\,} of coefficients in {\it growing\,} 
sets of partials of $A,\,B,\,C$ would play decisive roles in freezing or not of 
the values of new incoming pairs of component functions of the i.s.'s. 
Linear algebra packages would eventually come in handy. 

\vskip1mm
\noindent{\small {\textsc{Institute of Mathematics, University of Warsaw}
\newline\textit{Banach \,str.\,\,2\,, 02-097 Warsaw, Poland}}\\
\n{\textsc{E-mail address: }\textit{mormul@mimuw.edu.pl}}}
\vskip4mm
\n{\small{\textsc{Universit\'e de Savoie Mont Blanc (LAMA)\\
\textit{73370 Le Bourget-du-Lac, France}}\\ 
\n{\textsc{E-mail address: }\textit{fernand.pelletier@univ-smb.fr}}}

\end{document}
